\documentclass[10pt]{article}
\usepackage{amsfonts,amsmath,amsthm,amssymb,tikz}
\usepackage{graphics, epsfig}
\usepackage{color}
\usepackage{appendix}
\usepackage{verbatim}
\usepackage{ulem}
\usepackage{mathrsfs}
\usepackage[makeroom]{cancel}
\usepackage[margin=1.5in]{geometry}
\usepackage[colorlinks=true, urlcolor=blue, citecolor=red, linkcolor=blue, pdffitwindow=true, linktocpage, pdfpagelabels, bookmarksnumbered, bookmarksopen]{hyperref}
\usepackage{dsfont}
 \usepackage[usenames,dvipsnames]{pstricks}
 \usepackage{pst-grad} 
 \usepackage{pst-plot} 
\allowdisplaybreaks

\def\Xint#1{\mathchoice
{\XXint\displaystyle\textstyle{#1}}%
{\XXint\textstyle\scriptstyle{#1}}%
{\XXint\scriptstyle\scriptscriptstyle{#1}}%
{\XXint\scriptscriptstyle\scriptscriptstyle{#1}}%
\!\int}
\def\XXint#1#2#3{{\setbox0=\hbox{$#1{#2#3}{\int}$ }
\vcenter{\hbox{$#2#3$ }}\kern-.6\wd0}}

\def\dashint{\Xint-}

\newtheorem{proposition}{Proposition}[section]
\newtheorem{theorem}{Theorem}[section]
\newtheorem{lemma}{Lemma}[section]

\numberwithin{equation}{section}
\numberwithin{theorem}{section}
\numberwithin{proposition}{section}
\numberwithin{lemma}{section}
\numberwithin{remark}{section}
\newcommand{\osc}{\operatornamewithlimits{osc}}

\newcommand{\loc}{\operatorname{loc}}

\newcommand{\dist}{\operatorname{dist}}

\def\R{\mathbb{R}}
\def\Q{\mathcal{Q}}

\def\N{\mathbb{N}}

\def\dist{\mathrm{dist}}

\begin{document}

\title{A note on the point-wise behaviour of bounded solutions for a non-standard elliptic operator}
\author{
Laura Baldelli \& Simone Ciani \& Igor I. Skrypnik \&  Vincenzo Vespri\\ \\ Universit\`a degli Studi di Firenze,\\ Dipartimento di Matematica e Informatica "Ulisse Dini",\\   $laura.baldelli@unifi.it$  \& $vincenzo.vespri@unifi.it$, \\ \\
Institute of Applied Mathematics and Mechanics,\\
National Academy of Sciences of Ukraine,\\
$iskrypnik@iamm.donbass.com$\\\\ 
Technical University of Darmstadt, Department of Mathematics,\\
$ciani@mathematik.tu$-$darmstadt.de$}

\date{}
\maketitle
\vskip .4 truecm
\begin{abstract} \noindent
In this brief note we discuss local H\"older continuity for solutions to anisotropic elliptic equations of the type
\begin{equation*} \label{prototype}
\sum_{i=1}^s \partial_{ii} u+ \sum_{i=s+1}^N \partial_i \bigg(A_i(x,u,\nabla u) \bigg)  =0,\quad x \in \Omega \subset \subset \R^N \quad \text{ for } \quad 1\leq s \leq N-1,
\end{equation*}
where each operator $A_i$ behaves directionally as the singular $p$-Laplacian, $1< p < 2$ and the supercritical condition $p+(N-s)(p-2)>0$ holds true. We show that the Harnack inequality can be proved without the continuity of solutions and that  in turn this implies H\"older continuity of solutions.

\begin{flushright}
\it{To celebrate Jerry Goldstein's 80th genethliac}
\end{flushright}

\vskip.2truecm
\noindent {{\bf MSC 2020:}} 
35J75, 35K92, 35B65.
\vskip.2truecm
\noindent {{\bf Key Words}}: 
Anisotropic $p$-Laplacian, Singular Parabolic Equations, H\"older continuity, Intrinsic Scaling, $L^{1}-L^{\infty}$ estimates.
\end{abstract}

\newpage 
\section*{Introduction}\label{S:intro}
Along the following work we show that it is possible to provide an Harnack inequality independently of the continuity of solutions, for bounded weak solutions to some operators whose prototype  
\begin{equation}\label{prototipo}
    \sum_{i=1}^s \partial_{ii} u + \sum_{i=s+1}^N
 \partial_i \bigg( |\partial_i u|^{p-2}\partial_i u\bigg) =0, \quad \text{weakly in}\quad \Omega\subset \subset \R^N, \quad 1<p<2, 
 \end{equation} \noindent 
 has a {\it non-degenerate} behavior along the first $s$-variables, and a {\it singular} behavior on the last ones, with $s\in \mathbb{N}$, $1\le s\le N-1$. Moreover we show that H\"older continuity of solutions is a consequence of the sole Harnack inequality, following the heritage of Moser's approach (\cite{Moser61}, \cite{Moser}). The terms {\it degenerate} and {\it singular} are classic in evolutionary partial differential equations, for which we 
 refer for instance to the texts \cite{DB}, \cite{Vaz} and the references therein. 
 The type of operators in \eqref{prototipo} is widely studied to describe the steady states of non-Newtonian fluid with different directional diffusions, see for instance \cite{anto}, \cite{Lions}.
 On the other hand these adjectives are somehow of a minor use in the context of elliptic equations (see \cite{DB-book} for instance). In what follows we adopt this terminology, because it embodies our point of view on the study of regularity for these equations: we regard them as if they were parabolic, being the driving term the {\it nondegenerate} one. \vskip0.2cm \noindent 
 Before entering the details of the equation under study, we believe it is worth of interest to present the general problem, which represents an open mathematical challenge still after half a century. \vskip0.2cm \noindent 
Indeed, the anisotropic elliptic partial differential equation \eqref{prototipo} is a particular case of a more general kind of operators
\begin{equation} \label{general-anisotropic}
-\sum_{i=1}^N \partial_i A_i(x,u, \nabla u )= B(x,u,\nabla u), \quad \text{weakly in} \quad \Omega \subset \subset \R^N,    
\end{equation} \noindent 
where the maps $A_i,B: \Omega \times \R \times \R^N \rightarrow \R$ are Caratheodory functions (\cite{Giusti}, Definition 4.3) satisfying an anisotropic growth condition as 
 \begin{equation}\label{so fucking general}
     \begin{cases}
      \sum_{i=1}^N A_i(x,s,\xi) \cdot \xi_i \ge C_1 \sum_{i=1}^N |\xi_i|^{p_i} -C, \quad \text{for} \quad \xi \in \R^N, \\  \\
      |A_i(x,s,\xi)| \leq C_2 \sum_{i=1}^N |\xi_i|^{p_i-1}+C, \quad \text{for} \quad i\in \{ s+1,..,N\},\\ \\
      |B(x,s,\xi)| \leq C \sum_{i=1}^N |\xi_i|^{p_i-1}+C 
     \end{cases}
 \end{equation} \noindent where $C_1,C_2 >0$, $C \ge 0$ are given constants that we will always refer to as {\it the data}. The equation is called {\it homogeneous} when $C=0$. When $1<p_i\equiv p<\infty$ for all $i=1, \dots N$, the equation \eqref{general-anisotropic} reduces to an equation of $p$-growth, whose theory of regularity is reasonably complete (see for instance \cite{DB-book} Chap. X, \cite{Giusti}, \cite{LU}, \cite{Lions}, \cite{Moser61}, \cite{Trudy}). It does not coincide with the usual $p$-Laplacian, but with a similar equation called {\it pseudo}-$p$-Laplacian: the theory of 0th-order regulairty (as H\"older continuity and expansion of positivity) is similar to the $p$-Laplacian, but Lipschitz estimates become more difficult because the set of degeneration of derivatives is unboundend.  When $p_i$'s are different, the study of regularity was first studied in \cite{UU}, \cite{Marcellini1989}, \cite{Marcellini1991}, \cite{Marcellini1993}. See the recent surveys \cite{Marcellini-survey}-\cite{Mingione-Radulescu} for a more exhaustive list of references. A first condition to start a study of local regularity of solution is boundedness. It has been shown in \cite{FS}, \cite{Kolodii}, \cite{Korolev},
that local boundedness is inherent in the notion of local weak solution,  provided the following restriction of $p_i$'s, which is $p_i<N\bar{p}/(N-\bar{p})=:\bar{p}^*$ for all $i=1,\dots ,N$,
being $\bar{p}=(\sum_{i=1}^N 1/p_i)^{-1}$ the harmonic mean (see the subsection \ref{notations}). When this condition is violated, one can construct unbounded solutions, as local minimizers of energy integrals whose Euler-Lagrange equation fits into \eqref{general-anisotropic}-\eqref{so fucking general}, as 
 \[
 \mathcal{F}(u;\Omega)= \sum_{i=1}^N \frac{1}{p_i}\int_{\Omega} | \partial_i u |^{p_i} \, dx, \qquad \Omega \subset \subset \R^N, 
 \] with $p_i\equiv2$ for $i=1, \dots, N-1$ and $p_N$ big enough to violate the above condition $p_N<\bar{p}^*$. This is the content of \cite{Giaq-counter} and \cite{Marcellini-counter}. See also \cite{CMM} for a discussion under limit growth conditions. Lipschitz continuity has been proved by a double-iterative argument and some very particular choices of test functions in \cite{Brasco}. The bibliography is far from being complete.
 \vskip0.2cm  
 \noindent Regarding the H\"older continuity of solutions, the most general result insofar is a result of stability  \cite{DGV-remarks} (see also \cite{DICA}, \cite{DuzMarVes}, \cite{Eleuteri-Marcellini-Mascolo}, \cite{Eurica} for some interesting discussions around the subject). A Harnack inequality (and therefore, as we will show later, local continuity itself) is provided to the prototype elliptic equation in \cite{CMV}, but only for a range of exponents that is parabolic. See also the interesting work \cite{FVV} in the case of fast diffusion. The case of an equation similar to \eqref{prototipo} has been studied in \cite{LS}, for the {\it degenerate} case (i.e. $p>2$) and in \cite{LSV} for a similar structure than ours, but with the restriction $s=1$. We refer to \cite{CSV} for a discussion on the subject. In these two works the authors study the local properties of weak solutions to \eqref{prototipo} from a parabolic point of view. The main difference relies on the technique: while in the former the authors use a transformation which is strongly linked to condition $s=1$, in the latter the authors use the special structure of the equation, mainly that $p_i=2$ for $i=1,\dots, s$, with a particular choice of test functions that allows to write the {\it nondegenerate} part of the energy in a boundary-term. This boundary term is finally used to derive specific logarithmic estimates, as in \cite{DB-Chen}, \cite{DGV-Annali}. Both techniques spread new light on the interpretation of elliptic anisotropic operators, and as such can be seen as a starting point for a whole new theory. In this paper we mainly prove Harnack inequality, whose statement is given in Section \ref{def-main}, without assuming the property of continuity of solutions. Then, we get H\"older continuity in Section \ref{Hcont} by using the Harnack estimate, following the idea first developed by Moser in \cite{Moser} and then by DiBenedetto et al. in \cite{DGV-acta} for degenerate parabolic equations with rough coefficients. \vskip0.1cm 
\noindent The paper is structured as follows. In Section \ref{def-main} we introduce the functional setting, the definition of variational local weak solutions to equation \eqref{E1} and we state our main results. The description of the main analytic devices that are used along our work, included the consistency of the definition of solution with truncations are contained in Section \ref{Preliminaries}. Section \ref{semi} contains the point-wise behaviour of local weak solutions, while in Section \ref{Harnack-section}  we prove the Harnack inequality and finally in Section \ref{Hcont} we develop the H\"older continuity.  
    
\subsection*{Notations} \label{notations}
\begin{itemize}
\small{
\item[-]
If $\Omega$ is a measurable subset of $\R^{N}$, we denote by $|\Omega|$ its Lebesgue measure. We will write $\Omega \subset \subset \R^N$ when $\Omega$ is an open bounded set.
\item[-]
For $r>0$ and $\bar{x}=(\bar{x}',\bar{x}'') \in \R^{s} \times \R^{N-s}$, we denote by $B_{r}(\bar x)$ the ball of radius $r$ and center $\bar{x}$; the standard polydisc is denoted by $Q_{\theta,\rho}=B_{\theta}(\bar{x}') \times B_{\rho}(\bar{x}'') \subset \R^N$. Furthermore, by $w_s=|B_1(0')|$ and $w_{N-s}=|B_1(0'')|$ we denote the measures of the respective unit balls.
\item[-]
The symbol $\forall_{\text{ae}}$ stands for {\it -for almost every-} .   
\item[-]
For a measurable function $u:\Omega \to \R$, by $\inf u$ and $\sup u$ we understand the essential infimum and supremum, respectively. When $a\in \R$, we omit the domain when considering sub/super level sets, i.e. $\big[u\gtreqless a\big]=\big\{x\in \Omega: u(x)\gtreqless a\big\}$. We let  $\partial_i u$ denote the partial weak derivatives.
\item[-] We make the usual convention that a constant $\gamma>0$ depending only on the data, i.e. $\gamma=\gamma(N, 2, p, C_1, C_2, C)$, may vary from line to line along calculations.

\item[-] We consider positive numbers $(p_1,p_2, \dots, p_N)$ ordered increasingly and we call $\bar{p}= ( \sum_{i=1}^N 1/p_i)^{-1}$ and $\bar{p}^*= N \bar{p}/(N-\bar{p})$, respectively the {\it harmonic mean} of $p_i$'s and its {\it Sobolev's exponent}. Finally, we set \[ \lambda_{l} :=N(\bar{p}-2)+l\bar{p}, \qquad \text{and} \qquad \chi= \lambda_1= p + (N-s)(p-2).\]}

\end{itemize}

 \section{Definitions and Main Results} \label{def-main}
 \noindent 
Let $\Omega \subset \subset \R^N$ be an open bounded set with $N \ge 2$, and for $1<p<2$ and $1\leq s\leq N-1$ let us consider the elliptic partial differential equation
\begin{equation}\label{E1}
    \sum_{i=1}^s \partial_{ii} u + \sum_{i=s+1}^N
 \partial_i A_i(x,u, \nabla u)=0, \quad \qquad  \text{weakly in}\quad \Omega, \end{equation} where the Caratheodory\footnote{Measurable in $(u, \xi)$ for all $x \in \Omega$ and continuous in $x$ for a.e. $(u,\xi)\in \R\times \R^N$.} functions $A_i:\Omega\times  \R\times \R^N \rightarrow \R$ satisfy the following structure conditions for almost every $x \in \Omega$,
 \begin{equation}\label{E1-structure}
     \begin{cases}
      \sum_{i=s+1}^N A_i(x,u,\xi) \cdot \xi_i \ge C_1 \sum_{i=s+1}^N |\xi_i|^p -C, \quad \text{for} \quad \xi \in \R^N,\,\,   u \in \R,\,\,   x \in \Omega \\  \\
      |A_i(x,u,\xi)| \leq C_2 |\xi_i|^{p-1}+C, \quad \text{for} \quad i\in \{ s+1,..,N\},
     \end{cases}
 \end{equation} \noindent where $C_1,C_2 >0$, $C \ge 0$ are given constants that we will always refer to as {\it the data}. A function $u \in L^{\infty}_{loc}(\Omega) \cap W^{1,[2,p]}_{loc}(\Omega)$, where we define
\begin{equation*}\label{anisotropic-Sobolev-spaces}
\begin{aligned}
  &W_{loc}^{1,[2,p]}(\Omega):= \bigg{\{}u\in L^1_{loc}(\Omega)\,\,  | \,\,  \partial_i u \in L^2_{loc}(\Omega) \,\, \forall i=1,..,s, \quad \partial_i u \in L^p_{loc}(\Omega) \,\, \, \forall i=s+1,..,N   \bigg{\}}, \\
  &W^{1,[2,p]}_o(\Omega):= W^{1,1}_o(\Omega) \cap   W_{loc}^{1,[2,p]}(\Omega), 
  \end{aligned}
\end{equation*}\noindent
is called a local weak solution to \eqref{E1}-\eqref{E1-structure} if for each compact set $K \subset \subset \Omega$ it holds the following integral estimate
\begin{equation}\label{def-solution}
    \int \int_K \sum_{i=1}^s \partial_i u \, \partial_i \varphi \, dx + \int \int_K \sum_{s+1}^N A_i(x,u,\nabla u)\, \partial_i \varphi \, dx=0,\quad \forall \varphi \in W^{1,[2,p]}_o (K).
\end{equation}\noindent
\noindent All along the present work we will suppose that truncations of local weak solutions $u$ to \eqref{E1}-\eqref{E1-structure}, defined by 
\[
\pm (u-k)_{\pm} = \pm \max\{ \pm (u - k),\, \,  0\}, \quad \text{for} \quad  k \in \R,\]
preserve the property of being sub(super)-solutions. That is, for any $k \in \R$, every compact subset $K \subset\subset \Omega$, and $\psi \in W_o^{1,[2,p]}(K)$ we have
\begin{equation}\label{mon-energy}
\begin{aligned}
\int \int_{K} \bigg{\{} \sum_{i=1}^s  \partial_i (u-k)_{\pm} \partial_i\psi \, + \, \sum_{i=s+1}^N A_i(x,k\pm (u-k)_{\pm},\partial_i (u-k)_{\pm} ) \partial_i \psi \, \bigg{\}} dx \leq(\geq) 0.
\end{aligned}
\end{equation} \noindent 
Assumption \eqref{mon-energy} above is quite natural, see for instance Section \ref{Preliminaries}.\vskip0.2cm \noindent Let us fix some geometrical notations and convention. For a point $x_o\in \Omega$, let us denote it by  $x_o=(x_o',x_o'')$ where $x_o'\in \R^s$ and $x_o''\in \R^{N-s}$. Let $\theta,\rho>0$ be two parameters, and define the polydisc 
\[ Q_{\theta,\rho}(x_o):= B_{\theta}(x_o') \times B_{\rho}(x_o'').\]
 \noindent We will say $Q_{\theta,\rho}$ is an {\it intrinsic} polydisc when $\theta$ depends on the solution $u$ itself. We will call first $s$ variables the {\it nondegenerate} variables and last $(N-s)$ ones {\it singular} variables. Using this geometry we state our main result, an intrinsic form of Harnack's inequality.
 Let $\bar{p}$ be the harmonic mean of $p_i$s, i.e.
 \[\frac{1}{\bar{p}}= \frac{1}{N} \sum_{i=1}^N \frac{1}{p_i}=\frac{1}{N} \bigg(\frac{s}{2}+ \frac{N-s}{p}  \bigg), \quad \Rightarrow \quad \bar{p}=\frac{2Np}{2(N-s)+ps}. \]
The central point of the present work is that the following theorem can be obtained without the continuity of solutions.

 \begin{theorem}[Harnack estimate]\label{harnackTHM}
  Let $u$ be a bounded, nonnegative, local weak solution to \eqref{E1}-\eqref{E1-structure}-\eqref{mon-energy}. Let $x_o \in \Omega$ be a point such that $u(x_o)>0$ and $\rho>0$ small enough to allow the inclusion \begin{equation} \label{Mcal} Q_{\mathcal{M},\rho}(x_o) \subseteq \Omega, \quad \quad \text{being} \quad \mathcal{M}= ||u||_{L^{\infty}(\Omega)}^{(2-p)/2}\rho^{\frac{p}{2}}.\end{equation} \noindent
  Assume also that
  \begin{equation}\label{chi}
      \chi:=p+(N-s)(p-2)>0.
  \end{equation}
  
  \noindent Then there exist constants $K>1,\bar{\delta}_o\in (0,1)$ depending only on the data such that either
 \[  u(x_o) \leq K \rho,\] \noindent or
  \begin{equation}\label{Harnack}
      u(x_o) \leq K \inf_{Q_{\theta,\rho}(x_o)} u\,, \qquad \quad  \text{with} \quad \theta= \bar{\delta}_o u(x_o)^{\frac{2-p}{2}} \rho^{\frac{p}{2}}.
  \end{equation}
 \end{theorem} 
 
 \noindent The condition \eqref{Mcal} is part of the hypothesis and it ensures that polydisc $Q_{\theta,\rho}$ is contained in the set of definition of the equation $\Omega$. Estimate \eqref{Harnack} is valid for every $\rho>0$ small enough for the infimum to be taken in a defined set, and it is purely elliptic, i.e. it does not require a waiting "time". The intrinsic Harnack estimate \eqref{Harnack} above can be used, following an approach pioneered by J. Moser in \cite{Moser} and developed by many others (see for instance \cite{DGV-mono}), to prove that local weak solutions are locally H\"older continuous, with an estimate as the one below.
 \newpage

\begin{theorem}[H\"older Continuity] \label{HC}
Let $u$ be a bounded local weak solution to \eqref{E1} with structure conditions \eqref{E1-structure}-\eqref{mon-energy}. Moreover, let $p,s$ satisfy \eqref{chi}. Then $u$ is locally H\"older continuous. More precisely, there exist constants $\gamma>1$, $\alpha \in (0,1)$ depending only on the data $\{N,p,s,C_i, \, i=0,1,2\}$ such that for any compact set $K \subset \subset\Omega$ we have 
\begin{equation}\label{HC-estimate}
    |u(x)-u(y)| \leq \gamma \|u\|_{\infty,\Omega} \bigg( \frac{|x'-y'|^{\frac{2}{p}}\|u\|_{\infty,\Omega}^{\frac{p-2}{p}}+ |x''-y''|}{(2,p)-\dist(K, \partial \Omega)} \bigg)^{\alpha},  
\end{equation}\noindent for $x,y \in K$ and being 

\begin{equation}\label{p-dist}
    (2,p)-\dist(K,\partial \Omega)= \inf \bigg{\{}|x'-y'|^{\frac{2}{p}} \|u\|_{\infty,\Omega}^{\frac{p-2}{p}} + \, |x''-y''|,\, \, x\in K,\, y \in \partial \Omega \bigg{\}}.
\end{equation}
\end{theorem}

\section{Preliminaries and Tools of the Trade} \label{Preliminaries}
    
The Banach spaces of functions $W_{loc}^{1,[2,p]}(\Omega), W^{1,[2,p]}_o(\Omega)$ enjoy various properties of embedding and approximation by smooth functions, when their topology is induced by the norms
\[||u||_{W^{1,[2,p]}_{loc}(\Omega)}= ||u||_{L^{1}(\Omega)}+ \sum_{i=1}^s ||\partial_i u ||_{L^2(\Omega)}+\sum_{i=s+1}^N ||\partial_i u||_{L^p(\Omega)},\] 
and
\[||u||_{W^{1,[2,p]}_o(\Omega)}= \sum_{i=1}^s ||\partial_i u ||_{L^2(\Omega)}+\sum_{i=s+1}^N ||\partial_i u||_{L^p(\Omega)}.\]
    
\noindent For some references in this direction, we mention the following partial list \cite{Cianchi, KK, Schmeiser, Troisi, Tro2}. Notice the interesting approach of interpolation of \cite{Adams} to the anisotropic case. \vskip0.2cm \noindent Anisotropic spaces $W_{loc}^{1,[2,p]}(\Omega)$ are the spontaneous sets where to define variational solutions, even if the potential and viscosity approaches would not require a priori this integrability, at the price of a suitabe comparison principle. This is why our assumption \eqref{mon-energy} that truncations are subsolutions themselves is very natural. Let us consider for instance the monotone case, i.e. we assume that the following monotonicity property holds
  \begin{equation}\label{monotonicity}
      \sum_{i=s+1}^N \langle A_i(x,s,\eta)-A_i(x,s,\zeta), (\eta_i-\zeta_i) \rangle >0,
  \end{equation} \noindent  for $\eta, \zeta \in \R^{N-s}$ with $\eta_i \ne \zeta_i$ for all $i=s+1,..,N$ and $(x,s) \in \Omega \times \R$. The assumption \eqref{monotonicity} is natural when considering theory of existence and uniqueness of solutions to \eqref{prototipo} (see for instance \cite{Lions}, \cite{Showalter}), and eventually leads to the  assumption \eqref{mon-energy}.
\begin{proposition}
Let $u$ be a local weak solution to \eqref{E1} with structure conditions \eqref{E1-structure}. If monotonicity assumption \eqref{monotonicity} holds with 
\begin{equation} \label{integrability-Ass}
\partial_i A_i (x, u, 0) \in L^1_{loc}(\Omega) \qquad \text{for each} \quad i \in \{s+1, \dots,N\},
\end{equation} then property \eqref{mon-energy} holds.
\end{proposition}

\begin{proof}
Let $\varepsilon >0$ and test the equation \eqref{def-solution} with the test function $\varphi= \frac{\pm (u-k)_{\pm}}{\pm(u-k)_{\pm} + \varepsilon} \psi$, being $0 \leq \varphi, \psi \in W^{1,[2,p]}_o(\Omega)$ admissible, with $\varphi \leq \psi$ almost everywhere in $\Omega$. We prove the assertion for truncations from below $(u-k)_+$, the other case being similar. We observe that almost everywhere on the set $[u>k]$ we have 
\[ A_i (x,u,\nabla u) = A_i(x, k + (u-k)_{+}, \nabla (u-k)_{+}). \]



\noindent Applying the dominated convergence theorem we obtain
\begin{equation} \label{star-wars}
    \begin{aligned}
    \sum_{i=1}^s\int \int_K & \partial_i (u-k)_{+}\, \partial_i \psi\, dx + \sum_{i=s+1}^N\int \int_{K }  A_i(x, k + (u-k)_{+}, \nabla (u-k)_{+}) \, \partial_i \psi dx \, \\
    &\leq- \sum_{i=s+1}^N \liminf_{\varepsilon \downarrow 0} \varepsilon \int \int_K \frac{A_i(x, k + (u-k)_{+}, \nabla (u-k)_{+}) \, \partial_i (u-k)_+}{[(u-k)_+ +\varepsilon ]^2} \psi\, dx 
    \end{aligned}
\end{equation} \noindent
Now for each $i \in \{s+1,\dots,N\}$ we split the last term on the right of \eqref{star-wars} into the following decomposition:
\begin{equation*}
    \begin{aligned}
    \varepsilon  \int \int_K& \frac{A_i(x, k + (u-k)_{+}, \nabla (u-k)_{+}) \, \partial_i (u-k)_+}{[(u-k)_+ +\varepsilon ]^2}\, \psi \, dx\\
    &= \varepsilon \int \int_{K}  \frac{\langle A_i(x, k + (u-k)_{+}, \nabla (u-k)_{+})- A_i(x,k+(u-k)_+,0),\, \partial_i (u-k)_+\rangle }{[(u-k)_+ +\varepsilon ]^2}\, \psi\, dx \\
    & - \int \int_K \bigg{\{} \partial_i A_i(x, k+(u-k)_+, 0) \psi + A_i(x, k+(u-k)_+, 0) \partial_i \psi \bigg{\}} \,  \frac{(u-k)_+}{(u-k)_+ + \varepsilon} \, dx,     \end{aligned}
\end{equation*} \noindent and last term on the right is zero on the limit $\varepsilon \downarrow 0$, thanks to property \eqref{integrability-Ass}, divergence theorem and $\psi=0$ on $\partial K$. This means by monotonicity that the last term on the right of \eqref{star-wars} is negative and can be discarded, obtaining that $(u-k)_+$ is a sub-solution to an equation similar to \eqref{E1}. Being the structure conditions dependent mainly on the derivatives  $|\partial_i u| \chi_{[u>k]}= |\partial_i (u-k)_+|$, they remain unvaried.
\end{proof}

\noindent Here below we introduce some particular tools of the trade, which are particular properties of solutions to \eqref{mon-energy}- \eqref{E1-structure} that permitted in \cite{CSV} the achievement of Theorem \ref{harnackTHM}. The first result is about local boundedness of solutions, and its particular use should be deeply analyzed, as its condition on the exponent is the origin of the parabolic range \eqref{chi}. Observe that $\theta, \rho >0$ are parameters that can be chosen freely.
\begin{lemma} \label{boundedness}
Le $u$ be a locally bounded local weak solution to \eqref{E1}-\eqref{E1-structure}. Let $ 1\leq l \leq 2$ and 
\[ \lambda_{l} :=N(\bar{p}-2)+l\bar{p}>0,\] \noindent 
Then there exist constants $\gamma,C>0$ depending only on the data, such that for all polydiscs $Q_{2\theta,2\rho}\subset \Omega$ we have either
\[
\bigg( \frac{\theta^2}{\rho^p} \bigg)^{\frac{1}{2-p}}\leq C \rho, \quad \text{or}
\]
\[
    \sup_{Q_{\theta/2,\rho/2}} u  \leq \gamma  \bigg( \frac{\rho^p}{\theta^2} \bigg)^{( \frac{N-s}{p}) (\frac{\bar{p}}{ \lambda_{l} })}  \bigg( \dashint  \dashint_{Q_{\theta,\rho}} u_+^l\, dx  \bigg)^{\frac{\bar{p}}{ \lambda_{l} }} +\gamma \bigg(\frac{\theta^2}{\rho^p} \bigg)^{\frac{1}{2-p}}.\]\noindent 
\end{lemma}

\noindent Another fundamental tool for our analysis of local regularity is the following integral estimate, which can be seen as an Harnack estimate within the $L^1-L^{\infty}$ topology, and is typical of singular parabolic equations (see for instance \cite{DB}, Prop. 4.1 Chap VII).
 \begin{theorem}\label{l1-linftyTHM}
 Let $u$ be a nonnegative, bounded, local weak solution to \eqref{E1}-\eqref{E1-structure}. Fix a point $\bar{x} \in \Omega$ and numbers $\theta,\rho>0$ such that $Q_{8\theta,8\rho}(\bar{x}) \subset \Omega$. Then there exists a positive constant $\gamma$ depending only on the data such that either
 \begin{equation} \label{either}
     \bigg( \frac{\theta^2}{\rho^p}\bigg)^{\frac{1}{2-p}} \leq \rho,
 \end{equation}\noindent or
 \[
     \theta^{-s} \rho^{s-N} \int \int_{Q_{\theta,\rho}(\bar{x})} u \, dx \leq \gamma \bigg{\{}\inf_{B_{\frac{\theta}{2}}(\bar{x}')} \rho^{s-N} \bigg( \int_{B_{2\rho}(\bar{x}'')} u(\cdot, x'')\, dx'' \bigg)^{\frac{p}{\chi}}+ \bigg(\frac{\theta^2}{\rho^p} \bigg)^{\frac{1}{2-p}}     \bigg{\}}.
\] \noindent In this case, the exponent of the first integral on the right can be negative. If additionally property \eqref{chi} holds, then either we have \eqref{either} or
 \begin{equation} \label{l1-linfty}
     \sup_{Q_{\frac{\theta}{2},\frac{\rho}{2}}(\bar{x})} u \leq \gamma \bigg{\{} \bigg( \frac{\rho^p}{\theta^2} \bigg)^{\frac{N-s}{\chi}} \inf_{B_{\frac{\theta}{2}}(\bar{x}')}  \rho^{s-N} \bigg( \int_{B_{2\rho}(\bar{x}'')} u(\cdot, x'') dx''  \bigg)^{\frac{p}{\chi}} +\bigg(  \frac{\theta^2}{\rho^p} \bigg)^{\frac{1}{2-p}}  \bigg{\}}.
 \end{equation}
   \end{theorem}
\noindent 
The main tool to achieve both Harnack inequality and H\"older continuity in \cite{CSV} is a particular expansion of positivity, along the singular variables and valid for each relative measure of positivity.

\begin{theorem} \label{expansion}
  Let $\bar{x} \in \Omega$ and let $u\ge 0$ be a bounded, local weak solution to \eqref{E1}- \eqref{E1-structure}, satisfying \eqref{mon-energy}. Suppose that for a point $\bar{x} \in \Omega$ and numbers $M,\rho>0$ and $\nu \in (0,1)$ it holds
  \[
      |[u\leq M] \cap Q_{\theta,\rho}(\bar{x}) |\leq (1-\nu) |Q_{\theta,\rho}(\bar{x})|, \quad \text{for} \quad \theta= \rho^{\frac{p}{2}} (\delta M)^{\frac{2-p}{2}},
\] \noindent and $Q_{2\theta,2\rho}(\bar{x}) \subset \Omega$, for a number $\delta=\delta(\nu) \in (0,1)$. Then there exist constants $K>1$ and $\delta_o \in (0,1)$ depending only on the data and $\nu$ such that either
 \[
    M \leq K \rho,\] \noindent or
 \[
    u(x) \ge \delta_o M/2, \quad \forall_{ae}\,\, x \in Q_{\eta,2\rho}(\bar{x}), \quad \quad  \text{where} \quad \eta= (2\rho)^{\frac{p}{2}} (\delta_o M)^{\frac{2-p}{2}}. \]
 \end{theorem}

 \noindent Finally, we will use in a crucial way the following Lemma, proven in \cite{CMV}, following an approach that generalises a first idea on viscosity solutions (\cite{CC}). We believe that the viscosity solution approach would contribute significantly to the theory of anisotropic equations, because already the main problem of degeneracy on an infinite set can be seen in a method involving higher derivatives; see for instance \cite{dem} for the case $p_i \equiv p$ of the {\it pseudo} $p$-Laplace operator.

\begin{lemma}[Lemma 4.1 \cite{CMV}]\label{Krylov-Safonov}
Let $(X, {\rm d})$ be a quasi-metric space with quasi-metric constant $\gamma$ and $x_{0}\in X$. Let $\mathbb{B}_{\rho}(z)$ be the ball in $(X,d)$ of radius $\rho>0$ and center $z \in X$. Then for any $\beta>0$ there exists a constant $\omega=\omega(\gamma, \beta)>1$ such that for any  bounded function $u:\mathbb{B}_{1}(x_{0})\to \R$ with $u(x_{0})\ge 1$ there exist $ x\in \mathbb{B}_{1}(x_{0})$ and $r>0$ such that 
\[
\mathbb{B}_{r}(x)\subseteq \mathbb{B}_{1}(x_{0}),\qquad     r^{\beta}  \sup_{\mathbb{B}_{r}(x)} u\le \omega ,\qquad r^{\beta} u(x)\ge 1/\omega.\]
\end{lemma}

\section{Pointwise behaviour} \label{semi}
Lower semicontinuty of solutions to \eqref{E1}-\eqref{E1-structure} can be proven as a consequence of $L^{\infty}$ estimates (see for instance \cite{DMV},  \cite{Kuusi}) or more generally as a consequence of a theoretical maximum principle called De Giorgi-type Lemma (Lemma 2.4 in \cite{CSV}), in this respect we refer to \cite{CG, Liao}. Here we show an approach which is similar to \cite{DMV}, but it uses a different $L^{\infty}$ estimate, being $p<2$.

\begin{proposition} \label{semincontinuity}
Any bounded local weak solution $u$ of \eqref{E1} with \eqref{E1-structure} and $\bar{p}>2N/(N+1)$ in the supercritical range, has a lower semicontinuous representative.
\end{proposition}

\begin{proof}
We first consider $(\mathbb R^N, \mathscr{L}, d_M)$, where $\mathscr{L}$ is the Lebesgue's measure and $d_M$ is the distance defined as follows for any $0<M\in \mathbb{Q}$, $x,y\in \mathbb R^N$,
$$d_M(x,y)=max\{|x'-y'|^{2/p}M^{-2/p}, |x''-y''|\}.$$
Thus, notice that  $Q_{\rho^{p/2}M,\rho}(x_0)=\mathds{D}_\rho(x_0)$,
where $\mathds{D}_\rho(x_0)$ is the ball in $(\mathbb R^N, \mathscr{L}, d_M)$ of center $x_0\in \mathbb R^N$ and radius $\rho>0$. We observe that $(\mathbb R^N, \mathscr{L}, d_M)$, where $ \mathscr{L}$ is the standard Lebesgue measure, is a doubling metric measure space since 
\[ \mathscr{L}(\mathds{D}_{2\rho}(x_0))=|Q_{(2\rho)^{p/2}M,2\rho}(x_0)| \le C |Q_{\rho^{p/2}M,\rho}(x_0)|=\mathscr{L}(\mathds{D}_{\rho}(x_0)),
\] for an universal constant $C>0$. \newline
Now let $V_M(\Omega)$ the set of $L^1(\Omega)$-Lebesgue points of $(\Omega, \mathscr{L}, d_M)$, i.e. 
\[V_M(\Omega)= \bigg{\{} x \in \Omega: \, \, \dashint \dashint_{\mathds{D}_{\rho}(x)} |u(y)-u(x)|\, dy \rightarrow 0, \,\, \text{as}\,\, \rho \rightarrow 0   \bigg{\}}, \]
as (Theorem \cite{Hei}, Theorem 1.8]) guarantees $|V_M(\Omega)|=|\Omega|$ for each positive $M\in \mathbb{Q}$. Finally, we consider the full-measure set  
\[\bar{V}(\Omega)= \bigcap_{M\in \mathbb Q} V_M(\Omega),\qquad \qquad |\bar{V}_M(\Omega)|=|\Omega|.\]
Let us consider therefore a point $\bar{x}\in \bar{V}(\Omega)$, and let us modify $u \in W^{1,[2,p]}_{\loc}(\Omega)$ to be the representative in such a class. We need to divide the proof into two cases, due to non-homogeneity of the equation, by using Lemma \ref{boundedness}. \vskip0.2cm \noindent 
If $( \theta^2/\rho^p)^{1/(2-p)}\leq C \rho$, then by choosing $\theta= 2\|u\|_{\infty, \Omega}$ we obtain the condition 
\[
2\|u\|_{\infty, \Omega} \leq C \rho, \quad \Rightarrow \quad \osc_{Q_{\theta, \rho}(\bar{x})} u \leq C \rho,
\] and this reduction of oscillation implies Lipschitz continuity of solutions (see Section \ref{Hcont}). \newline
If on the other hand $( \theta^2/\rho^p)^{1/(2-p)}> C \rho$, then the following estimate is valid by choosing $l=1$ and $\theta=\rho^{p/2}M$, where $M>0$ to be chosen, so that we get 
\begin{equation}\label{2.12}
\sup_{Q_{\theta/2,\rho/2}(\bar{x})} u \leq \gamma M^{-2(N-s)\bar{p}/p\lambda_1} \bigg( \dashint \dashint_{Q_{\theta,\rho}(\bar{x})} u(x)_+ dx \bigg)^{\bar{p}/\lambda_1} +\gamma M^{2/(2-p)},
\end{equation}
where  $\lambda_1>0$ since $\bar{p}>2N/(N+1)$ is in force.

\noindent Now, being $\bar{x}\in \bar{V}(\Omega)$, for any $\varepsilon_1>0$, there exists $\rho=\rho(\varepsilon_1)>0$ such that
\begin{equation}\label{lebx0}
\dashint \dashint_{Q_{\theta,\rho}(\bar{x})} |u(x)-u(\bar{x})| dx <\varepsilon_1.
\end{equation}
Since $u-u(\bar{x})$ is a solution of an equation of the type of \eqref{E1}, then for every $\varepsilon>0$, by using \eqref{lebx0} in \eqref{2.12}, we have
\begin{equation}\label{ast}
\sup_{Q_{\theta/2,\rho/2}(\bar{x})} (u(\bar{x})-u)_+ \leq \gamma \biggl(M^{-2(N-s)\bar{p}/p\lambda_1} \varepsilon_1^{\bar{p}/\lambda_1} + M^{2/(2-p)}\biggr)\le \frac{\varepsilon}{4},\end{equation}
by choosing appropriately
\begin{equation}\label{M}
\varepsilon_1=\biggl(\frac{\varepsilon M^{2(N-s)\bar{p}/p\lambda_1}}{2\gamma}\biggr)^{\lambda_1/\bar{p}}, \quad M=(\varepsilon/2\gamma)^{(2-p)/2}.
\end{equation}

\noindent Now we claim the lower semicontinuity, that is the following
$$u(\bar{x})\le ess \liminf_{x\to\bar{x}} u(x)= ess \lim_{\rho\to0}\inf_{Q_{\theta, \rho}(\bar{x})} u(x).$$
We proceed by contradiction, supposing that there exists $\bar{\varepsilon}>0$ and $R_0>0$ such that for all $0<\rho<R_0$, we have
$$u(\bar{x})-\lim_{\rho\to0}\inf_{Q_{\theta/2, \rho/2}(\bar{x})} u(x)=\bar{\varepsilon}>0.$$
Then, choosing again $\theta=\rho^{p/2}M$ and $M$ as in $\eqref{M}_2$, so that, if $\varepsilon=\bar{\varepsilon}$ defined before, we reach the required contradiction since, by \eqref{ast}
$$\bar{\varepsilon}=\sup_{Q_{\frac{\theta}{2},\frac{\rho}{2}}(\bar{x})} (u(\bar{x})-u)_+\le \frac{\bar{\varepsilon}}{4}.$$
The proof is so concluded.
\end{proof}

\section{Proof of Harnack inequality \ref{Harnack} without continuity} \label{Harnack-section}
Consider a bounded local weak solution $u$ to \eqref{E1} in $\Omega$, let $x_o \in \Omega$ be a point where\footnote{Notice that this point-wise value can be identified thanks to Proposition \ref{semincontinuity}.} $u(x_o)>0$, and for any $\rho \ge 0$ appropriate to satisfy \eqref{Mcal}, construct the cylinder $Q_{\theta,\rho} (x_o)$ such that
\[Q_{\theta,\rho} (x_o)= B_{\theta}(x_o') \times B_{\rho}(x_o'')\subset \subset Q_{\mathcal{M}, \rho}(x_o) \subseteq \Omega, \qquad \text{with} \qquad \theta= u(x_o)^{\frac{2-p}{2}}\rho^{\frac{p}{2}},\]
with $\mathcal{M}$ defined in \eqref{Mcal}.
Now we proceed to normalize the solution $u: Q_{\theta,\rho}(x_o)\rightarrow \R$ with the transformation 
\begin{equation} \label{trans}
    v(x',x'')= \frac{1}{u(x_o)} u\bigg(\frac{x'-x_o'}{\theta}, \frac{x''-x_o''}{\rho}\bigg),
\end{equation} \noindent 
that is a bounded local weak solution to an equation similar to \eqref{E1} in $Q_{1,1}$, i.e. solves 
\begin{equation}\label{E-trans}
    \sum_{i=1}^s \partial_{ii} v(x)+ \sum_{i=s+1}^N \partial_i \tilde{A}_i(x,v,\nabla v)=0, \qquad \text{weakly in} \qquad Q_{1,1}= B_1(0') \times B_1(0''),
\end{equation} \noindent with structure conditions
\begin{equation}\label{E-trans-structure}
\begin{cases}
\sum_{i=s+1}^N \tilde{A}_i(x,s,\xi) \xi_i \ge \tilde{C_1} \sum_{i=s+1}^N |\xi_i|^p- \tilde{C}^p,\\
|\tilde{A}_i(x,s,\xi)|\leq \tilde{C_2} \sum_{i=s+1}^N |\partial_i v|^{p-1} + \tilde{C}^{p-1},
\end{cases} \text{being} \qquad \tilde{C}=(C \rho)/ u(x_o).
\end{equation} \noindent For parameters $\lambda \in (0,1)$ and $\beta>1$ to be chosen later, we consider the equation
\begin{equation} \label{Krylov-Safonov-EQ}
    \sup_{x''\in B_{\lambda_0}} u(0',x'') = \frac{ (1-\lambda_0)^{-\beta}}{2}=:M/2.
\end{equation} \noindent Then, by continuity, the proof in \cite{CSV} exhibited a maximal root $\lambda_0$ of the equation \eqref{Krylov-Safonov-EQ} together with a point $\bar{x}''\in B_{\lambda_0}(0'')$ such that, by defining suitable $r=r(\lambda_0)$, the authors, exploiting a popular argument of Krylov and Safonov (see \cite{KS}), obtain the estimate  \begin{equation}\label{l1}
\frac{M}{2}\leq v(0',\bar{x}'') \leq \sup_{B_{r}(\bar{x}'')} v(0', \, \cdot)\leq \gamma M.
\end{equation} From this estimate on, the proof in \cite{CSV} follows the lines of \cite{DGV-Annali} and uses twice \eqref{l1-linfty} to arrive at a measure estimate of the kind 
\begin{equation}\label{measure-harnack}
    |[v(x',\cdot) \leq \varepsilon M] \cap B_{r/2}(\bar{x}'')| \leq (1-\alpha) |B_{r/2}|, \qquad \text{for} \quad \varepsilon>0 \quad \text{and} \quad \alpha \in (0,1),
\end{equation} \noindent and finally use the expansion of positivity along the singular variables, i.e. Theorem \ref{expansion}, to obtain a lower bound 
\[ v>\tilde{M}, \quad \text{in} \quad Q_{\tilde{\eta},1}\, , \quad \text{being} \quad \tilde{\eta}= \delta_o(\tilde{M}) r^{p/2}, \]
for some constants $\tilde{M},\delta_o(\tilde{M})\in (0,1)$ depending only on $\beta$ and the data. At this point the proof is concluded by transforming back the function $v$ into $u$, and consequently the estimate above transforms into the final claim:
\[u(x_o) \leq \tilde{M}^{-1} u(x), \qquad x \in Q_{\eta,\rho}(x_o), \qquad \text{with} \quad \eta= \delta_o u(x_o)^{\frac{2-p}{2}} \rho^{\frac{p}{2}}. \]
\vskip0.2cm \noindent 
Now, the argument that follows shows that we can as well obtain estimate above  \eqref{l1} without using the continuity of $u$. We use indeed Lemma \ref{Krylov-Safonov} for the bounded function $v(0',\cdot):B_1(0'')\to \R$ of the sole $x''$-variable, in the euclidean metric space $(B_1(0''),d_{e})$, where $d_e(x'',y'')= \sqrt{\sum_{i=s+1}^N (x_i''-y_i'')^2}$. Therefore for each $\beta >0$ we find $\bar{x}'' \in {B}_1(0'')$ and ${\gamma},r>0$ such that by setting $M=r^{-\beta}$ we obtain
\[\frac{M}{\gamma} \leq v (0',\bar{x}'') \leq \sup_{B_r(\bar{x}'')} v(0', \, \cdot) \leq \gamma M.  \]
But this is precisely \eqref{l1} by redefining the constant $M$, and the proof can now carried on the same way. 
\qed

\section{H\"older Continuity} \label{Hcont}
Let $K$ be a compact subset of $\Omega$, consider $(2,p)-\dist(K,\partial \Omega)$ as in \eqref{p-dist}.
For any point $y_o \in K$ we construct the polydisc $Q_{\theta_o, R}(y_o)$ with
\begin{equation}\label{Rdef}
R= [(2,p)-\dist (K, \partial \Omega)]/2, \qquad \theta_o=R^{p/2} \omega_o^{\frac{2-p}{2}},\qquad \omega_o= 2 \|u\|_{\infty, \Omega}. \end{equation}
We claim that $Q_{\theta_o, R}(y_o)\subset\Omega$
 by the choice of $R$ and $\theta_o$ above. Indeed, for $ z \in Q_{\theta_o, R}(y_o)$, it holds
\[|z'-y_o'| \leq \theta_o=R^{p/2} \omega_o^{\frac{2-p}{2}}  \leq \inf \{|x'-y'|, \, \, x \in K, \, y \in \partial \Omega\},\]and
\[|z''-y_o''|\le R \leq \inf \{|x''-y''|,\, \, x \in K, \, y \in \partial \Omega  \}/2.\]
Thus, by \eqref{p-dist}, we get $Q_{\theta_o, R}(y_o)\subset \Omega$. \newline \noindent 
Next, we consider a generic point $x_o  \in K$, and we show that we can reduce the proof of H\"older Continuity to the assumption
\begin{equation} \label{accommodation} 
x_o \in  Q_{\theta_o, R}(y_o).
\end{equation}
Indeed, if this is not the case, namely $x_o \notin Q_{\theta_o, R}(y_o)$, we obtain Lipschitz continuity. In particular, assume either
\[|x_o'-y_o'|>\theta_o=R^{p/2} \omega_o^{\frac{2-p}{2}}, \quad \text{or} \quad |x_o''-y_o''|>R,\]
then, by computing straightforward, we obtain the Lipschitz condition
\[|u(x_o)-u(y_o)| \leq 2 \omega_o\leq 4 \omega_o \bigg( \frac{|x_o'-y_o'|^{\frac{2}{p}}\omega_o^{\frac{p-2}{p}}+ |x_o''-y_o''|}{R}  \bigg),\]
where we used the definition of $\omega_o$, $R$.
So we proceed by assuming \eqref{accommodation} and, before proving the H\"older Continuity, we show that there exists a sequence of intrinsic polydiscs $\Q_{\theta_n,\rho_n}$ with center in $y_o$ where the oscillation of $u$ can be controlled uniformly.

\begin{proposition}[Oscillation Decay] \label{osc-decay}
Let $u$ be a bounded local weak solution to \eqref{E1} with structure conditions \eqref{E1-structure}-\eqref{mon-energy}. Assume also $y_o\in K$ and define $\omega_o$ and $R$ as in \eqref{Rdef}. There exist a constant $\delta \in (0,1)$ such that if we define
\[
\begin{cases}
\rho_n=\delta^n R,\\
\omega_n= \delta^n \omega_o,
\end{cases}    
\begin{cases}
\delta= 4K/(4K+1)\in (0,1),\\
\theta_n= \bar{\delta_o} \rho_n^{p/2} \omega_n^{\frac{2-p}{2}}=2^{\frac{p-2}{2}}\bar{\delta_o} \delta^{n} R^{p/2} \omega_o^{\frac{2-p}{2}},
\end{cases} \]
where $K>1$ defined in Theorem \ref{harnackTHM},
then 
\[\Q_{n+1}\subset \Q_n, \quad \text{for} \quad \Q_n= y_o+ \Q_{\theta_n, \rho_n}, \quad \Q_0=y_o+ \Q_{\theta_o,R},\] \noindent 
and 
\begin{equation} \label{osc}\osc_{\Q_n} u \leq \delta^n \omega_o.\end{equation} 
\end{proposition}

\begin{proof}   
We proceed by induction observing first that the first step is achieved from the definition of $\omega_o$ so that we have $ \osc_{\Q_0} u \leq \omega_o$. Moreover, we observe that the polydiscs $\Q_n$ for $n\in\mathbb N$ are monotonic. Indeed, since $\omega_1 \leq \omega_0$ and $\theta_1 \leq \theta_o$ then $\Q_1 \subset \Q_0$ and this is sufficient to establish $\Q_{n+1}\subset \Q_n$ for all $n \in \N$. So we assume that \eqref{osc} is valid until step $n$ and we prove it for the $n+1$-th. We also assume, by contradiction, that $\osc_{\Q_{n+1}}u >\omega_{n+1}$. Now we define \[M_n=\sup_{\Q_n} u, \quad m_n=\inf_{\Q_n} u,\]
and we suppose that one of the two following estimates holds
\begin{equation}\label{dicot}M_n- u(y_o) > \omega_{n+1}/4, \qquad \text{or} \qquad u(y_o)-m_n> \omega_{n+1}/4.\end{equation}
This can be assumed because otherwise we have
\[\osc_{\Q_{n+1}}u \leq \osc_{\Q_n} u \leq    \omega_{n+1},\] 
leading to an absurd. Let us suppose that \eqref{dicot}$_2$ is valid, the other case is similar. Let \eqref{Mcal} and \eqref{chi} be valid, then we apply Theorem \ref{harnackTHM}, so that \eqref{Harnack} is valid for the non-negative function $u-m_n$, which is still a solution of a similar equation to \eqref{E1} with \eqref{E1-structure}. So we obtain 
\begin{equation} \label{contr}(4K)^{-1} \omega_{n+1} < K^{-1}( u(y_o)-m_{n} ) \leq \inf_{Q_{\theta,R}(y_o)} (u -m_n)\leq \inf_{\Q_{n+1}} (u - m_n)=0,\end{equation} \noindent provided that $\Q_{n+1} \subset Q_{\theta, R}( y_o)$, i.e.
\[ \bar{\delta_o}  \rho_{n+1}^{p/2} \omega_{n+1}^{\frac{2-p}{2}}\leq \bar{\delta_o} (u(y_o)-m_n)^{\frac{2-p}{2}} R^{p/2}, \quad \text{and} \quad \delta^{n+1} R \leq R. \]
The first inequality is satisfied exactly thanks to our assumption $\omega_{n+1}<2(u(y_o)-m_n)$. 
Finally, since $\Q_{n+1}\subset \Q_n$, by \eqref{contr}, the inductive hypothesis and definition of $\delta$ and $\omega_n$, the absurd inequality is reached
\[\osc_{\Q_{n+1}} u < M_n-m_n-(4K)^{-1}\omega_{n+1} \leq
\omega_n -(4K)^{-1}\omega_{n+1} =
\omega_{n+1} \bigg(\frac{1}{\delta}-\frac{1}{4K} \bigg)= \omega_{n+1}.\]




\end{proof}

\noindent Now we are ready to prove the H\"older Continuity, namely Theorem \ref{HC}, whose statement is given in Section \ref{def-main}.

\subsubsection*{Conclusion of the Proof of Theorem \ref{HC}}
For $x_o \in\Q_{\theta_o,R}(y_o)$ let $n \in \mathbb{N}$ be the last number such that $x_o \in \Q_n$ but $x_o \not\in \Q_{n+1}$. The latter implies that either
\[|x_o'-y_o'| >\theta_{n+1}=  \gamma_o \delta^{n} R^{p/2} \omega_o^{\frac{2-p}{2}} ,\quad \text{that is} \quad  \delta^{n2/p}<\gamma\, \,  \frac{|x_o'-y_o'|^{2/p} \omega_o^{\frac{p-2}{p}}}{R}\]
or
\[|x_o''-y_o''| >\rho_{n+1}= \delta^{n} R,\quad \text{that is} \quad \delta^{n}<\frac{|x_o''-y_o''|}{R} \] 
for a constant $\gamma_o= \gamma_o(\delta_o)>0$ depending only on the data. Being $p<2$ and estimating quantities that are smaller than $1$, we have 
\[\delta^{2n/p} \leq \gamma  \bigg(\frac{|x_o''-y_o''|}{R}+ \frac{|x_o'-y_o'|^{2/p} \omega_o^{\frac{p-2}{p}}}{R} \bigg) ,\]
\[\qquad \qquad \Rightarrow \delta^n \leq \gamma \bigg( \frac{|x_o''-y_o''|+|x_o'-y_o'|^{2/p} \omega_o^{\frac{p-2}{p}}}{R} \bigg)^{p/2}.  \] This, together with condition $ x_o \in \Q_n$ and Proposition \ref{osc-decay} gives 
\[|u(x_o)-u(y_o)| \leq \delta^n \omega_o \leq \gamma \omega_o  \bigg(  \frac{|x_o''-y_o''|+|x_o'-y_o'|^{2/p} \omega_o^{\frac{p-2}{p}}}{R} \bigg)^{\frac{p}{2}},\]
the proof is concluded.
\qed

\section*{Acknowledgements}
L. Baldelli and V. Vespri are supported by GNAMPA group of INdAM. S. Ciani is supported from the department of mathematics of Technical University of Darmstadt.

\end{document}